\documentclass [12pt]{article}

\usepackage{amsfonts,amsmath,amssymb}
\usepackage{geometry}
\newcommand{\li}{\mathop{\rm li}\nolimits}
\newcommand{\Li}{\mathop{\rm Li}\nolimits}
\def\qed{\hfill\nobreak\Box}
\newcommand{\Res}{\mathop{\rm Res}\limits}
\def\Z{{\mathbb Z}}

\begin{document}
\begin{center}
\large{
\textbf{Integrals Over Polytopes,}
\textbf{Multiple Zeta Values and Polylogarithms, and Euler's Constant}
}
\end{center}

\vskip 5mm

\begin{center}
Jonathan Sondow \\
209 West 97th Street \\
New York City \\
New York 10025 USA \\
jsondow(at)alumni.princeton.edu \\

\vskip 2mm

and

\vskip 2mm

Sergey Zlobin \\
Faculty of Mechanics and Mathematics \\
Moscow State University \\
Leninskie Gory, Moscow 119899 RUSSIA \\
sirg\_zlobin(at)mail.ru
\end{center}

\vskip 5mm

\begin{center}
\textbf{Abstract}
\end{center}

Let $T$ be the triangle with vertices (1,0), (0,1), (1,1). We study certain
integrals over $T$, one of which was computed by Euler. We give expressions for
them both as a linear combination of multiple zeta values, and as a
polynomial in single zeta values. We obtain asymptotic expansions of the
integrals, and of sums of certain multiple zeta values with constant weight.
We also give related expressions for Euler's constant. In the final section,
we evaluate more general integrals~-- one is a Chen (Drinfeld-Kontsevich)
iterated integral~-- over some polytopes that are higher-dimensional analogs of $T$. This
leads to a relation between certain multiple polylogarithm values and
multiple zeta values.

\tableofcontents

\section{Introduction}

Let $T$ be the triangle defined by
\[
T: = {\left\{ {{\left. {(x,y) \in [0,1]^{2}} \right|}x + y \ge 1}
\right\}},
\]
with vertices (1,0), (0,1), (1,1). In this paper we study the integral over
$T$
\begin{equation}
\label{eq1}
I_{n} : = {\int\!\!\!\int_{T} {{\frac{{( - \ln xy)^{n}}}{{xy}}}dxdy}}
\end{equation}
for $n = - 1,0,1,2,\dots$. We also consider integrals over several
polytopes that are higher-dimensional analogs of $T$.

Euler computed an iterated integral equivalent to $I_{0} $, and found that
\[
I_{0} = {\int\!\!\!\int_{T} {{\frac{{dxdy}}{{xy}}}}}  = {\int_{0}^{1}
{{\frac{{1}}{{x}}}{\int_{1 - x}^{1} {{\frac{{dy}}{{y}}}dx}}} }  =
{\int_{0}^{1} {{\frac{{ - \ln (1 - x)}}{{x}}}dx}}  = {\int_{0}^{1}
{{\sum\limits_{r = 1}^{\infty}  {{\frac{{x^{r - 1}}}{{r}}}dx}}} }  =
{\sum\limits_{r = 1}^{\infty}  {{\frac{{1}}{{r^{2}}}}}}  = \zeta (2).
\]
Using integration by parts he derived formula (\ref{eq19}), and used it to calculate
$\zeta (2)$ correctly to six decimals---see [\textbf{5}, Section 1.2],
[\textbf{7}, pp. 43-45].

We generalize Euler's result to $n = 0, 1, 2, \dots$ by showing that
$I_{n}$ is equal to an integer linear combination of 
\textit{multiple zeta values}

\[
\zeta (s_{1}, \dots, s_{l} ): = {\sum\limits_{n_{1} > n_{2} > \cdots
> n_{l} > 0} {{\frac{{1}}{{n_{1}^{s_{1}}  \cdots
n_{l}^{s_{l}}} } }}}
\]

\noindent
of \textit{weight} $s_{1} + \cdots + s_{l} = n + 2$. We also express $I_{n} $ as
a polynomial in single zeta values.

\bigskip

\noindent
\textbf{Theorem 1} \textit{Let} $n \ge 0$ \textit{ be an integer.} \\
(i) \textit{Then}
\begin{equation}
\label{eq2}
I_{n} = n! {\sum\limits_{k = 0}^{n} {\zeta (n - k + 2,\{1\}_{k}
)}} ,
\end{equation}
\textit{where} $\{1\}_{k}$ \textit{denotes}
$1, 1, \dots, 1$ (\textit{k times}). \\
(ii) \textit{Moreover,} $I_{n} $ \textit{is equal to an explicit polynomial
of several variables with rational coefficients in the values
of the Riemann zeta function}
$\zeta (2), \zeta (3), \dots, \zeta (n + 2).$

\bigskip

\noindent
(Theorem 1, Corollary 1, and Lemma 1 were obtained by the second author in
[\textbf{17}].) The proof is given in Section 2, along with the explicit
formula. Examples are
$$
I_{0} = \zeta (2),
\quad
I_{1} = \zeta (3) + \zeta (2,1) = 2\zeta (3),
\quad
I_{2} = 2\left( {\zeta (4) + \zeta (3,1) + \zeta (2,1,1)} \right) =
{\frac{{9}}{{2}}}\zeta (4),
$$
\begin{equation}
\label{missing1}
I_{3} = 6\left( {\zeta (5) + \zeta (4,1) + \zeta (3,1,1) + \zeta (2,1,1,1)}
\right) = 36\zeta (5) - 12\zeta (2)\zeta (3).
\end{equation}

The cases $n = 0, 1,$ and 2 are particularly simple.

\bigskip

\noindent
\textbf{Corollary 1} \textit{For} $n = 0, 1$, \textit{and} 2,
\textit{the integral} $I_{n} $ \textit{is a rational multiple of} $\zeta (n + 2)$.

\bigskip

For $n = 0$ and 1, this also follows from Beukers' [\textbf{2}] formulas
for $\zeta (2)$ and $\zeta (3)$ as integrals over the unit square
\[
S: = [0,1]^{2}.
\]
Namely, the change of variables $x = X, y = 1 - XY$ transforms both $I_{0}
$ into

\[
I_{0} = {\int\!\!\!\int_{T} {{\frac{{dxdy}}{{x y}}}}}  =
{\int\!\!\!\int_{S} {{\frac{{dXdY}}{{1 - XY}}}}}  = \zeta (2)
\]

\noindent
and ${\frac{{1}}{{2}}}I_{1} $ into
\begin{align*}
{\frac{{1}}{{2}}}I_{1} = {\frac{{1}}{{2}}}{\int\!\!\!\int_{T} {{\frac{{ -
\ln xy}}{{x y}}}}} dxdy & = {\int\!\!\!\int_{T} {{\frac{{ - \ln
x}}{{x y}}}}} dxdy \\
& = {\int\!\!\!\int_{S} {\;{\frac{{ - \ln X}}{{1
- XY}}}}} dXdY = {\frac{{1}}{{2}}}{\int\!\!\!\int_{S} {\;{\frac{{ - \ln
XY}}{{1 - XY}}}}} dXdY = \zeta (3).
\end{align*}

Here is an outline of the proof of Theorem 1. We first prove

\bigskip

\noindent
\textbf{Lemma 1} \textit{If} $k \ge 0$\textit{ and} $l \ge 0$
\textit{ are integers, then}
\begin{equation}
\label{eq3}
I_{k,l} : = {\int\!\!\!\int_{T} {{\frac{{( - \ln x)^{k}( - \ln
y)^{l}}}{{xy}}}}} dxdy = k!l!\zeta (l + 2,\{1\}_{k} ).
\end{equation}
\textit{If in addition} $l \ge 1$\textit{, then}

\begin{equation}
\label{eq4}
J_{k,l} : = {\int_{0}^{1} {{\frac{{( - \ln (1 - x))^{k}}}{{1 - x}}}( - \ln
x)^{l}dx}}  = k!l!\zeta (l + 1,\{1\}_{k} ).
\end{equation}

\noindent
Expanding $( - \ln x - \ln y)^{n}$, part (i) follows immediately. To prove
the lemma, we show that $(l + 1)I_{k,l} = J_{k,l + 1} $, and then evaluate
the integral $J_{k,l} $. Part (ii) of Theorem 1 follows, using a formula in
[\textbf{9}] for $J_{k,l} $ in terms of single zeta values.

As an application, we obtain an explicit version of a result in
[\textbf{3}].

\bigskip

\noindent
\textbf{Corollary 2} \textit{If} $n \ge 2$ \textit{and} $k \ge 0$\textit{,
then the multiple zeta value} $\zeta (n,\{1\}_{k} )$
\textit{ can be explicitly represented as a polynomial
of several variables with rational coefficients in the single zeta values}
$\zeta (2), \zeta (3), \dots, \zeta (n + k).$

\bigskip

Lemma 1 also affords a simple proof of a special case of the duality theorem
for multiple zeta values (see, for example, [\textbf{5}, Section 2.8]).

\bigskip

\noindent
\textbf{Corollary 3} \textit{If} $k \ge 0$\textit{ and} $l \ge 0$\textit{, then} $\zeta (k + 2,\{1\}_{l} ) = \zeta (l
+ 2,\{1\}_{k} ).$

\bigskip

For instance, $\zeta (2,1) = \zeta (3)$ and $\zeta (2,1,1) = \zeta (4)$.
Using these equalities, we give a second proof of Corollary 1.  However, unlike the cases $n = 0$ and 1, we do not have a proof of the
case $n = 2$ of Corollary 1 that does not use Theorem 1.

On the basis of numerical evidence and examples such as (\ref{missing1}),
we make the

\bigskip

\noindent
\textbf{Conjecture 1} \textit{The integral} $I_{n} $ \textit{is} not\textit{ a rational multiple of} $\zeta (n + 2)$\textit{ when} $n > 2$.

\bigskip

This has not been proved for a single value of $n$. However, using Theorem 1
(ii), we give a conditional proof for all $n = 3, 4, \dots$, assuming a
standard conjecture (see, for example, [\textbf{16}, Introduction]).

\bigskip

\noindent
\textbf{Theorem 2} \textit{If the numbers}
$\pi, \zeta (3), \zeta (5), \zeta (7), \zeta (9), \dots$
\textit{ are algebraically independent over the rationals, then Conjecture}
1\textit{ is true.}

\bigskip

Using a lemma we prove which gives an asymptotic expansion for the
coefficients of the Taylor series of certain meromorphic functions (Lemma
2), we estimate $I_{n} $ for $n$ large.

\bigskip

\noindent
\textbf{Theorem 3} \textit{The asymptotic equivalence}
\begin{equation}
\label{missing2}
I_{n} \sim 2n! \quad (n \to \infty )
\end{equation}

\noindent
\textit{holds. More precisely, the following asymptotic expansion is valid}:

$$
{\frac{{I_{n}}} {{n!}}} \approx 2 + {\frac{{6}}{{2^{n + 2}}}} +
{\frac{{20}}{{3^{n + 2}}}} + {\frac{{70}}{{4^{n + 2}}}} + \cdots \quad
(n \to \infty ),
$$

\noindent
\textit{where the numerator of the k-th term is}
$\binom{2k}{k}$, \textit{for} $k = 1, 2, \dots$.

\bigskip

This in turn gives an estimate for the sum of the multiple zeta values
$\zeta (m - k,\{1\}_{k} )$ of constant weight $m$.

\bigskip

\noindent
\textbf{Corollary 4}\textit{ The average of the multiple zeta values}
$\zeta (m), \zeta (m - 1,1), \dots, \zeta(2,\{1\}_{m - 2} )$
\textit{ is asymptotic to} ${{2} \mathord{\left/ {\vphantom {{2} {m}}} \right.
\kern-\nulldelimiterspace} {m}}$\textit{ as m tends to infinity}.
\textit{ In fact, the following asymptotic expansion holds}:

$$
{\sum\limits_{k = 0}^{m - 2} {\zeta (m - k,\{1\}_{k} )}}  \approx 2 +
{\frac{{6}}{{2^{m}}}} + {\frac{{20}}{{3^{m}}}} + {\frac{{70}}{{4^{m}}}} +
\cdots \quad (m \to \infty ).
$$

Another application of Theorem 3 is a curious result.

\bigskip

\noindent
\textbf{Corollary 5} \textit{The series}

\[
{\sum\limits_{n = 0}^{\infty}  {( - 1)^{n}{\frac{{I_{n}}} {{n!}}}}}
\]

\noindent
\textit{diverges, but is Abel summable to} ${{1} \mathord{\left/ {\vphantom {{1} {2}}} \right.
\kern-\nulldelimiterspace} {2}}.$

\bigskip

Let us now go \verb|"|down\verb|"| from $I_{0} $ to $I_{ - 1} $.

\bigskip

\noindent
\textbf{Question} \textit{Can one evaluate the integral}

\begin{equation}
\label{eq5}
I_{-1} = {\int\!\!\!\int_{T} {{\frac{{dxdy}}{{xy( - \ln xy)}}}}}  =
1.7330025 \dots
\end{equation}

\noindent
\textit{in terms of more familiar constants?}

\bigskip

Surprisingly, it turns out that $I_{ - 1} $ involves all the integrals
$I_{0}, I_{1}, I_{2}, \dots$ (hence all multiple zeta values $\zeta
(m,\{1\}_{k} )$ for $m \ge 2$ and $k \ge 0)$.

\bigskip

\noindent
\textbf{Theorem 4} \textit{If} $\li$
\textit{is the logarithmic integral function, then}

\[
I_{ - 1} = {\sum\limits_{n = 0}^{\infty}  {( - 1)^{n}{\frac{{I_{n}}} {{(n +
1)!}}}}}  + {\int_{0}^{1} {{\frac{{\li(x - x^{2})}}{{x}}}dx}}  + 1.
\]

\noindent
(Compare the convergent series here with the divergent series in Corollary
5.)

We now transform the double integral $I_{ - 1} $ into single integrals, one
involving the \textit{generalized binomial coefficient}

\[
\binom{s}{t} : = {\frac{{\Gamma (s + 1)}}{{\Gamma (t + 1)\Gamma (s
- t + 1)}}}.
\]

\noindent
\textbf{Proposition 1} \textit{The following integral formulas for} $I_{ - 1} $ \textit{are valid}:

\begin{equation}
\label{eq6}
I_{ - 1} = {\int_{0}^{\infty}  {\left( {1 - {\frac{{1}}{\binom
{2t}{t}}}} \right){\frac{{dt}}{{t^{2}}}}}}  = {\int_{0}^{1}
{\ln \left( {1 + {\frac{{\ln (1 - x)}}{{\ln x}}}} \right){\frac{{dx}}{{x}}}}
}.
\end{equation}

Expanding the first integrand in a power series, we find that the $n$th
coefficient involves the integral $I_{n} $.

\bigskip

\noindent
\textbf{Theorem 5} \textit{If} $0 < {\left| t \right|} < 1$
\textit{, then}

\begin{equation}
\label{eq7}
\left( 1 - \frac{{1}}{\binom{2t}{t}} \right){\frac{{1}}{{t^{2}}}} =
{\sum\limits_{n = 0}^{\infty}  {( - 1)^{n}{\frac{{I_{n}}} {{n!}}}t^{n}}} .
\end{equation}

\noindent
An application is Corollary 5.

We now relate $I_{ - 1} $ to \textit{Euler's constant} $\gamma $, which is defined as the limit

\[
\gamma = {\mathop {\lim} \limits_{n \to \infty}}  \left( {1 +
{\frac{{1}}{{2}}} + \cdots + {\frac{{1}}{{n}}} - \ln n}
\right).
\]

\noindent
If one thinks of $\gamma $ as \verb|"|$\zeta (1)$,\verb|"|
then from the formulas $I_{2} =
{\frac{{9}}{{2}}}\zeta (4)$, $I_{1} = 2\zeta (3)$, and $I_{0} = \zeta (2)$
one might expect that $I_{ - 1} $ involves $\gamma $. This is also suggested
by the similarity between the double integral (\ref{eq5}) for $I_{ - 1} $ and the
double integral for Euler's constant [\textbf{12}], [\textbf{14}]

\begin{equation}
\label{eq8}
\gamma = {\int\!\!\!\int_{S} {{\frac{{1 - X}}{{(1 - XY)( - \ln XY)}}}}
}dXdY.
\end{equation}

Formula (\ref{eq6}) leads to another, related connection between $I_{ - 1} $ and
$\gamma $. Namely, when $t = n$ is a positive integer, $\binom{2t}{t}$
is the central binomial coefficient $\binom{2n}{n}$,
which figures in the formulas for Euler's constant
$$
\binom{2n}{n} \gamma =
A_{n} - L_{n} + {\int\!\!\!\int_{S}
{{\frac{{\left( {X(1 - X)Y(1 - Y)} \right)^{n}}}{{(1 - XY)( - \ln XY)}}}}
}dXdY \quad (n \ge 1)
$$

\noindent
and

\[
\gamma = \frac{A_{n} - L_{n}} {\binom{2n}{n}} +
O\left( {{\frac{{1}}{{2^{6n}\sqrt {n}}} }}
\right) \quad (n \to \infty),
\]

\noindent
where $A_{n} $ is a certain rational number and $L_{n} $ is a particular
linear form in logarithms [\textbf{12}].

If in (\ref{eq8}) we perform the change of variables $X = x, Y = {{(1 - y)}
\mathord{\left/ {\vphantom {{(1 - y)} {x}}} \right.
\kern-\nulldelimiterspace} {x}}$, we obtain an integral over the triangle
$T$ for Euler's constant,

\begin{equation}
\label{eq9}
\gamma = {\int\!\!\!\int_{T} {{\frac{{1 - x}}{{xy( - \ln (1 - y))}}}}
}dxdy,
\end{equation}

\noindent
analogous to the triangle integral (\ref{eq5}) for $I_{ - 1} $.

We find an analog for $\gamma $ of the first integral for $I_{ - 1} $ in
(\ref{eq6}), which involves the generalized binomial coefficient
$\binom{2t}{t}$. (There exist classical analogs for $\gamma $ of the
second integral, which involves logarithms.)

\bigskip

\noindent
\textbf{Proposition 2} \textit{The following formula for Euler's constant is valid}:

\[
\gamma = {\int_{0}^{\infty}  {{\sum\limits_{k = 2}^{\infty}
{{\frac{{1}}{{k^{2} \binom{t + k}{k}}}}dt}}} } .
\]

\noindent
As an application, if we integrate termwise, and exponentiate the resulting
series, we recover \textit{Ser's infinite product for} $e^{\gamma} $ [\textbf{11}] (rediscovered in
[\textbf{13}], [\textbf{15}]):

\[
e^{\,\gamma}  = {\prod\limits_{k = 2}^{\infty}  {\left( {{\prod\limits_{j =
1}^{k} {j^{( - 1)^{j} \binom{k - 1}{j - 1}}}}}
\right)^{{{1\,} \mathord{\left/ {\vphantom
{{1\,} {\,k}}} \right. \kern-\nulldelimiterspace} {\,k}}}}}  = \left(
{{\frac{{2}}{{1}}}} \right)^{{{1\,} \mathord{\left/ {\vphantom {{1\,}
{\,2}}} \right. \kern-\nulldelimiterspace} {\,2}}}\left( {{\frac{{2^{2}}}{{1
\cdot 3}}}} \right)^{{{1\,} \mathord{\left/ {\vphantom {{1\,} {\,3}}}
\right. \kern-\nulldelimiterspace} {\,3}}}\left( {{\frac{{2^{3} \cdot 4}}{{1
\cdot 3^{{\kern 1pt} 3}}}}} \right)^{{{1\,} \mathord{\left/ {\vphantom
{{1\,} {\,4}}} \right. \kern-\nulldelimiterspace} {\,4}}}\left(
{{\frac{{2^{4} \cdot 4^{4}}}{{1 \cdot 3^{{\kern 1pt} 6} \cdot 5}}}}
\right)^{{{1\,} \mathord{\left/ {\vphantom {{1\,} {\,5}}} \right.
\kern-\nulldelimiterspace} {\,5}}} \cdots.
\]

The rest of the paper is organized as follows. In Sections 2 and 3 we
establish the non-asymptotic and asymptotic results, respectively, on $I_{n}
$ for $n \ge 0$. The applications to multiple zeta values are proved in
Section 4, and in Section 5 we prove the formulas for $I_{ - 1} $ and
$\gamma $. The final section is devoted to generalizing $I_{n} $ to
integrals over higher-dimensional analogs of the triangle $T$;
one is a Chen (Drinfeld-Kontsevich) iterated integral (see Remark 3). An application
is a relation between certain multiple polylogarithm values and multiple
zeta values (Corollary 7).

\section{The integral $I_{n}$ for $n \ge 0$}

We prove the non-asymptotic results on $I_{0}, I_{1}, I_{2}, \dots$
stated in the Introduction.

\bigskip

\noindent
\textbf{Lemma 1} \textit{If} $k \ge 0$\textit{ and} $l \ge 0$\textit{ are integers, then}

\[
I_{k,l} : = {\int\!\!\!\int_{T} {{\frac{{( - \ln x)^{k}( - \ln
y)^{l}}}{{xy}}}}} dxdy = k! l! \zeta (l + 2,\{1\}_{k} ).
\]

\noindent
\textit{If in addition} $l \ge 1$\textit{, then}

\[
J_{k,l} : = {\int_{0}^{1} {{\frac{{( - \ln (1 - x))^{k}}}{{1 - x}}}( - \ln
x)^{l}dx}}  = k! l! \zeta (l + 1,\{1\}_{k} ).
\]

\noindent
\textit{Proof.} We have

\[
I_{k,l} = {\int_{0}^{1} {{\frac{{( - \ln x)^{k}}}{{x}}}{\int_{1 - x}^{1}
{{\frac{{( - \ln y)^{l}}}{{y}}}dydx}}} }  = {\int_{0}^{1} {{\frac{{( - \ln
x)^{k}}}{{x}}} \cdot {\frac{{( - \ln (1 - x))^{l + 1}}}{{l + 1}}}dx}} {\rm
.}
\]

\noindent
Replacing $x$ with $1 - x$, we see that

\[
I_{k,l} = {\frac{{J_{k,l + 1}}} {{l + 1}}}.
\]

\noindent
Thus (\ref{eq3}) follows from (\ref{eq4}). To prove (\ref{eq4}), we multiply the formula
[\textbf{16}, Section 1]

\[
( - \ln (1 - x))^{k} = k!{\sum\limits_{n_{1} > n_{2} > \cdots >
n_{k} > 0} {{\frac{{x^{n_{1}}} }{{n_{1} \cdots n_{k}}} }}}
\]

\noindent
by $(1 - x)^{ - 1} = 1 + x + x^{2} + \cdots $, and substitute the
resulting series

\[
{\frac{{( - \ln (1 - x))^{k}}}{{1 - x}}} = k!{\sum\limits_{m \ge n_{1} >
n_{2} > \cdots > n_{k} > 0} {{\frac{{x^{m}}}{{n_{1} \cdot \cdot
\cdot n_{k}}} }}}
\]

\noindent
into the integral (\ref{eq4}) for $J_{k,l} $. We then integrate termwise, using the
fact that

\[
{\int_{0}^{1} {x^{m}( - \ln x)^{l}dx}}  = {\frac{{l!}}{{(m + 1)^{l +
1}}}}.
\]

\noindent
The result is

\[
J_{k,l} = k!l!{\sum\limits_{m > n_{1} > n_{2} > \cdots > n_{k}
> 0} {{\frac{{1}}{{m^{l + 1}n_{1} \cdots n_{k}}} }}}  =
k!l!\zeta (l + 1,\{1\}_{k} ),
\]

\noindent
and the lemma follows. $\qed$

\bigskip

\noindent
\textbf{Theorem 1} \textit{If} $n \ge 0$\textit{, then} $I_{n} $\textit{ can be expressed both}

\noindent
(i) \textit{in terms of} \textit{multiple zeta values as}

\begin{equation}
\label{eq10}
I_{n} = n! {\sum\limits_{k = 0}^{n} {\zeta (n - k + 2,\{1\}_{k}
)}}
\end{equation}
(ii) \textit{and in terms of single zeta values as}

\begin{equation}
\label{eq11}
I_{n} = \sum\limits_{k = 0}^{n} \binom{n}{k}
\frac{J_{k,n - k + 1}} {n - k + 1},
\end{equation}

\noindent
\textit{where the integral} $J_{k,n - k + 1} $,\textit{ defined in} (\ref{eq4}),\textit{ is given by the formula} [\textbf{9}]

\begin{equation}
\label{eq12}
J_{k,l} = k!l!{\sum\limits_{p = 1}^{l} {{\frac{{( - 1)^{p +
1}}}{{p!}}}{\sum\limits_{t_{i}}  {{\frac{{\zeta (t_{1} ) \cdots
\zeta (t_{p} )}}{{t_{1} \cdots t_{p}}} }{\sum\limits_{l_{i}}
\binom{t_{1}}{l_{1}} \cdots 
\binom{t_{p}}{l_{p}} }} }} } ,
\end{equation}

\noindent
\textit{the sum on} $t_{i} $\textit{ being taken over all sets of integers}
$\{t_{1}, \dots, t_{p} \}$\textit{ with}

\[
t_{i} > 1,\;\;\;{\sum\limits_{i = 1}^{p} {t_{i} = k + l + 1}} ,
\]

\noindent
\textit{and the sum on} $l_{i} $\textit{ over all sets of integers}
$\{l_{1}, \dots, l_{p} \}$ \textit{ with}

\[
0 < l_{i} < t_{i} ,\;\;\; {\sum\limits_{i = 1}^{p} {l_{i} = l}} .
\]

\noindent
\textit{Proof.} Expanding $( - \ln xy)^{n} = ( - \ln x - \ln y)^{n}$ in the definition (\ref{eq1})
of $I_{n} $, and applying (\ref{eq3}), gives (\ref{eq10}). Hence, using (\ref{eq4}), formula (\ref{eq11})
holds. Finally, the evaluation (\ref{eq12}) of the integral (\ref{eq4}) for $J_{k,l} $ is
proved in [\textbf{9}]. $\qed$

\bigskip

\noindent
\textbf{Corollary 1} \textit{For} $n = 0,\;1$,\textit{ and} 2,
\textit{the integral} $I_{n} $ \textit{is a rational multiple of}
$\zeta (n + 2)$.

\bigskip

We give two proofs. The first is short, but uses Theorem 1 (ii), whose proof
depends on [\textbf{9}]. The second is longer, but is self-contained (except
for a formula due to Euler): it uses Theorem 1 (i) and Corollary 3, whose
proofs do not rely on other papers.

\bigskip

\noindent
\textit{Proof} 1. For $n = 0,\;1$, and 2, formulas (\ref{eq11}) and (\ref{eq12}) yield $I_{0} = J_{0,1} =
\zeta (2)$ and $I_{1} = {\frac{{1}}{{2}}}J_{0,2} + J_{1,1} = 2\zeta (3)$ and
$I_{2} = {\frac{{1}}{{3}}}J_{0,3} + {\frac{{1}}{{2}}}J_{1,2} + J_{2,1} =
{\frac{{9}}{{4}}}\zeta (4)$. $\qed$

\bigskip

\noindent
\textit{Proof} 2. In the Introduction, we showed that $I_{0} = \zeta (2)$. Using the same
method, together with the formula $\smallint _{0}^{1} x^{k - 1}( - \ln
x)\,dx = k^{ - 2}$, we obtain

\[
I_{1} = 2{\int\!\!\!\int_{T} {{\frac{{ - \ln x}}{{xy}}}}} dxdy =
2{\int_{0}^{1} {{\frac{{\ln (1 - x)}}{{x}}}\ln x\,dx}}  = 2{\sum\limits_{k =
1}^{\infty}  {{\frac{{1}}{{k}}}{\int_{0}^{1} {x^{k - 1}( - \ln x)\,dx}}} }
= 2\zeta (3).
\]

Alternatively, $I_{0} = \zeta (2)$ and $I_{1} = \zeta (3) + \zeta (2,1)$ by
Theorem 1 (i), and $\zeta (2,1) = \zeta (3)$ by Corollary 3.

In order to prove that $I_{2} = {\frac{{9}}{{2}}}\zeta (4)$, it suffices, by
Theorem 1 (i) and Corollary 3, to apply Euler's formula $\zeta (3,1) =
{\frac{{1}}{{4}}}\zeta (4)$. (For the latter, take $n = 3$ in his equation
(9.5) of [\textbf{1}, p. 252].) $\qed$

\bigskip

\noindent
\textbf{Theorem 2} \textit{If the numbers}
$\pi, \zeta (3), \zeta (5), \zeta (7), \zeta (9), \dots$
\textit{ are algebraically independent over the rationals, then}
$I_{n} $ \textit{is not a rational multiple of} $\zeta (n + 2)$
\textit{ when} $n > 2.$

\bigskip

\noindent
\textit{Proof.} First take the case $n = 3m - 2$, with $m > 1$. The integral $I_{n} $ is
equal to a linear combination (\ref{eq11}) of integrals $J_{k,l} $ with positive
coefficients. Each $J_{k,l} $ is equal to a polynomial (\ref{eq12}) of several
variables in single zeta values. Now in (\ref{eq12}) the monomial $\zeta (3)^{m}$
appears only when $p = m$, and then its coefficient is nonzero and has sign
$( - 1)^{m + 1}$. Hence in the expression for $I_{n} $ the coefficient of
$\zeta (3)^{m}$ is nonzero. It follows, using the hypothesis, that $I_{n} $
cannot be a rational multiple of $\zeta (n + 2) = \zeta (3m)$.

The cases $n = 3m + 3$ and $n = 3m + 5$, with $m > 0$, are similar: consider
the monomials $\zeta (3)^{m}\zeta (5)$ and $\zeta (3)^{m} \zeta (7)$,
respectively. The remaining cases $n = 3$ and $n = 5$ can be handled by
direct calculation, completing the proof. $\qed$

\bigskip

\noindent
\textbf{Theorem 5} \textit{If} $0 < {\left| t \right|} < 1$\textit{, then}

\begin{equation}
\label{eq13}
{\frac{{1}}{{t^{2}}}}\left( {1 - {\frac{{1}}{\binom
{2t}{t}}}} \right) = {\sum\limits_{n = 0}^{\infty}  {( -
1)^{n}{\frac{{I_{n}}} {{n!}}}t^{n}}} .
\end{equation}

\noindent
\textit{Proof.} The generating function

\begin{equation}
\label{eq14}
{\sum\limits_{k,l \ge 0} {x^{k + 1}y^{l + 1}\zeta (l + 2,\{1\}_{k} )}}  = 1
- {\frac{{\Gamma (1 - x)\Gamma (1 - y)}}{{\Gamma (1 - x - y)}}}
\end{equation}

\noindent
(compare (\ref{eq16})) is derived in [\textbf{3}]. If $x = y = - t$, then the series
converges when ${\left| t \right|} < 1$. Setting
$k + l = n$, we obtain

\[
t^{2}{\sum\limits_{n = 0}^{\infty}  {( - 1)^{n}t^{n}{\sum\limits_{k = 0}^{n}
{\zeta (n - k + 2,\{1\}_{k} )}}} }  = 1 - {\frac{{1}}{\binom
{2t}{t}}}.
\]

\noindent
Applying (\ref{eq2}), the theorem follows. $\qed$

\section{Asymptotic expansion of $I_{n}$}

Using Theorem 5 and the next lemma, we estimate the integral $I_{n} $ when
$n$ is large.

\bigskip

\noindent
\textbf{Lemma 2}\textit{ Suppose that the function} $f(z)$\textit{ is meromorphic in the complex plane and has only simple poles}
$z_{1} ,\;z_{2} ,\dots$ \textit{, with residues} $r_{1} ,\;r_{2}
,\dots$\textit{, respectively. If}
$0 < |z_{1}| \le
|z_{2}| \le \cdots $\textit{, then the coefficients of the Taylor series}

\[
f(z) = {\sum\limits_{n = 0}^{\infty}  {a_{n} z^{n}}}
\]

\noindent
\textit{have the asymptotic expansion}

$$
a_{n} \approx - {\frac{{r_{1}}} {{z_{1}^{n + 1}}} } - {\frac{{r_{2}
}}{{z_{2}^{n + 1}}} } - \cdots \quad (n \to \infty ).
$$

\noindent
\textbf{Remark 1} Recall [\textbf{6}, Section 1.5] that the last formula
means that, for every fixed positive integer $k$,
\begin{equation}
\label{missing3}
a_{n} = - {\frac{{r_{1}}} {{z_{1}^{n + 1}}} } - \cdots -
{\frac{{r_{k}}} {{z_{k}^{n + 1}}} } + O\left( {{\frac{{r_{k + 1}}} {{z_{k +
1}^{n + 1}}} }} \right) \quad (n \to \infty ).
\end{equation}

\noindent
\textit{Proof of Lemma} 2$.$ A meromorphic function has only
finitely many poles in any bounded region,
so for each $k \ge 1$ there exists $l > k$ with
$|z_{l}| < |z_{l + 1}|$. Note that the only singularities of the function

\begin{equation}
\label{eq15}
f(z) - {\sum\limits_{j = 1}^{l} {{\frac{{r_{j}}} {{z - z_{j}}} }}}  =
{\sum\limits_{n = 0}^{\infty}  {\left( {a_{n} + {\frac{{r_{1}}} {{z_{1}^{n +
1}}} } + \cdots + {\frac{{r_{l}}} {{z_{l}^{n + 1}}} }}
\right)z^{n}}}
\end{equation}

\noindent
are $z_{l + 1} ,\;z_{l + 2}, \dots$. It follows, using the inequalities
$0 < |z_{1}| \le |z_{2}| \le \cdots $, that the radius
of convergence of the series (\ref{eq15}) is equal to
$|z_{l + 1}|$. As $|z_{l}| < |z_{l + 1}|$, we may substitute
$z = z_{l} $ into the series. Therefore,

\[
{\mathop {\lim} \limits_{n \to \infty}}  \left( {a_{n} + {\frac{{r_{1}
}}{{z_{1}^{n + 1}}} } + \cdots + {\frac{{r_{k}}} {{z_{k}^{n + 1}
}}} + {\frac{{r_{k + 1}}} {{z_{k + 1}^{n + 1}}} } + \cdots +
{\frac{{r_{l}}} {{z_{l}^{n + 1}}} }} \right)z_{l}^{n} = 0.
\]

\noindent
Since $|z_{k + 1}| \le |z_{k + 2}| \le \cdots \le
|z_{l}|$, the limit implies the
asymptotic formula (\ref{missing3}). This proves the lemma. $\qed$

\bigskip

\noindent
\textbf{Theorem 3} \textit{We have} $I_{n} \sim 2n!
$\textit{ as n tends to infinity. More precisely, the following asymptotic expansion is valid}:

$$
{\frac{{I_{n}}} {{n!}}} \approx 2 + {\frac{{6}}{{2^{n + 2}}}} +
{\frac{{20}}{{3^{n + 2}}}} + {\frac{{70}}{{4^{n + 2}}}} + \cdots \quad
(n \to \infty ),
$$

\noindent
\textit{where the numerator of the k-th term is}
$\binom{2k}{k}$ \textit{for} $k = 1, 2, \dots$.

\bigskip

\noindent
\textit{Proof.} Denote the function on the left side of (\ref{eq7}) by $f(t)$.
Aside from a
removable singularity at $t = 0$, the singularities of $f(t)$ are simple
poles at $t = - 1, - 2, \dots$. A calculation shows that the residue at
$t = - k$ is equal to

\[
\Res(f; - k) = \Res\left( {{\frac{{1}}{{t^{2}}}}\left( {1 - {\frac{{\Gamma (t
+ 1)^{2}}}{{\Gamma (2t + 1)}}}} \right); - k} \right) = -
{\frac{{1}}{{k^{2}}}}{\mathop {\lim} \limits_{t \to - k}} (t +
k){\frac{{\Gamma (t + 1)^{2}}}{{\Gamma (2t + 1)}}} = {\frac{{1}}{{k}}}
\binom{2k}{k}
\]

\noindent
for $k = 1,\;2,\dots$. Using Theorem 5 and Lemma 2, the second statement
(which implies the first) follows. $\qed$

\bigskip

\noindent
\textbf{Corollary 5} \textit{The series}

\[
{\sum\limits_{n = 0}^{\infty}  {( - 1)^{n}{\frac{{I_{n}}} {{n!}}}}}
\]

\noindent
\textit{diverges, but is Abel summable to}
${{1} \mathord{\left/ {\vphantom {{1} {2}}} \right.
\kern-\nulldelimiterspace} {2}}.$

\bigskip

\noindent
\textit{Proof.} The divergence follows from (\ref{missing2}).
Letting $t \to 1^{ -} $ in (\ref{eq13}), we obtain
the desired Abel summation. $\qed$

\section{Applications to multiple zeta values}

Using the results obtained on $I_{n} $, we study multiple zeta values of the
form $\zeta (m,\{1\}_{k} )$.

\bigskip

\noindent
\textbf{Corollary 2} \textit{If} $m \ge 2$ and $k \ge 0$, \textit{
then the multiple zeta value} $\zeta (m,\{1\}_{k} )
$\textit{ can be explicitly represented as a polynomial of several
variables with rational coefficients in the single zeta values}
$\zeta(2), \zeta (3), \dots, \zeta (m + k).$

\bigskip

\noindent
\textit{Proof.} Set $l = m - 1$ in (\ref{eq4}) and (\ref{eq12}). $\qed$

\bigskip

\noindent
\textbf{Remark 2} This result, including the polynomial formula (at least
implicitly), was first obtained in [\textbf{3}], using the generating
function (see [\textbf{3}] for the equivalence with (\ref{eq14}))

\begin{equation}
\label{eq16}
{\sum\limits_{k,l \ge 0} {x^{k + 1}y^{l + 1}\zeta (l + 2,\{1\}_{k} )}}  = 1
- \exp \left( {{\sum\limits_{n = 2}^{\infty}  {{\frac{{x^{n} + y^{n} - (x +
y)^{n}}}{{n}}}\zeta (n)}}}  \right).
\end{equation}

\noindent
\textbf{Corollary 3} \textit{If} $k \ge 0$\textit{ and} $l \ge 0$\textit{, then} $\zeta (k + 2,\{1\}_{l} ) = \zeta (l
+ 2,\{1\}_{k} ).$

\bigskip

\noindent
\textit{Proof.} Making the change of variables $x,y \to y,x$ in the integral (\ref{eq3}), the
symmetry of the triangle $T$ yields $I_{k,l} = I_{l,k} $. Using Lemma 1, the
result follows. $\qed$

\bigskip

\noindent
\textbf{Corollary 4}
\textit{ The average of the multiple zeta values} $\zeta (m),\;\zeta (m - 1,1),\dots\;,\;\zeta
(2,\{1\}_{m - 2} )$\textit{ is asymptotic to} ${{2} \mathord{\left/ {\vphantom {{2} {m}}} \right.
\kern-\nulldelimiterspace} {m}}$\textit{ as m tends to infinity}.\textit{ In fact, the following asymptotic expansion holds}:

$$
{\sum\limits_{k = 0}^{m - 2} {\zeta (m - k,\{1\}_{k} )}}  \approx 2 +
{\frac{{6}}{{2^{m}}}} + {\frac{{20}}{{3^{m}}}} + {\frac{{70}}{{4^{m}}}} +
\cdots  \quad (m \to \infty ).
$$

\noindent
\textit{Proof.} Setting $m = n + 2$ in Theorems 1 and 3 gives the desired expansion. It
follows that the average in question is asymptotic to ${\frac{{2}}{{m -
1}}}\sim {\frac{{2}}{{m}}}$ as $m$ tends to infinity. $\qed$

\section{The integral $I_{-1}$ and Euler's constant}

We prove the results on $I_{ - 1} $ and $\gamma$ stated in the Introduction.

\bigskip

\noindent
\textbf{Proposition 1} \textit{The following single integral formulas
for the double integral} $I_{ - 1} $ \textit{are valid}:

\[
I_{ - 1} = {\int_{0}^{\infty}  {\left( {1 - {\frac{{1}}{\binom{2t}{t}}}}
\right){\frac{{dt}}{{t^{2}}}}}}  = {\int_{0}^{1}
{\ln \left( {1 + {\frac{{\ln (1 - x)}}{{\ln x}}}} \right){\frac{{dx}}{{x}}}}
}.
\]

\noindent
\textit{Proof.} In (\ref{eq5}) make the substitution

\begin{equation}
\label{eq17}
 - {\frac{{1}}{{\ln xy}}} = {\int_{0}^{\infty}  {(xy)^{t}dt}}
\end{equation}

\noindent
and change the order of integration, obtaining

\begin{align*}
I_{ - 1} & = {\int_{0}^{\infty}  {{\int_{0}^{1} {{\int_{1 - x}^{1} {(xy)^{t -
1}dydx}}} } dt}}  = {\int_{0}^{\infty}  {{\int_{0}^{1} {\left( {x^{t - 1} -
x^{t - 1}(1 - x)^{t}} \right)\,}} dx{\frac{{dt}}{{t}}}}} \\ & =
{\int_{0}^{\infty}  {\left( {{\frac{{1}}{{t}}} - {\frac{{\Gamma (t)\Gamma (t
+ 1)}}{{\Gamma (2t + 1)}}}} \right)\,{\frac{{dt}}{{t}}}}} ,
\end{align*}

\noindent
using Euler's integral for the beta function. Replacing $\Gamma (t)$ with
$t^{ - 1}\Gamma (t + 1)$, the first equality follows. To see that the second
integral is also equal to $I_{ - 1} $, integrate with respect to $y$ in (\ref{eq5}).
$\qed$

\bigskip

\noindent
\textbf{Theorem 4} \textit{If} $li$\textit{ is the logarithmic integral, then}

\[
I_{ - 1} = {\sum\limits_{n = 0}^{\infty} {( - 1)^{n}{\frac{{I_{n}}} {{(n +
1)!}}}}}  + {\int_{0}^{1} {{\frac{{\li(x - x^{2})}}{{x}}}dx}}  + 1.
\]

\noindent
\textit{Proof.} Make the substitution

\[
 - {\frac{{1}}{{\ln xy}}} = {\int_{0}^{1} {(xy)^{t}dt}}  - {\frac{{xy}}{{\ln
xy}}}
\]

\noindent
in (\ref{eq5}). Using the proof of Proposition 1, we get

\[
I_{ - 1} = {\int_{0}^{1} {\left( {1 - {\frac{{1}}{
\binom{2t}{t}}}} \right){\frac{{dt}}{{t^{2}}}}}}  -
{\int\!\!\!\int_{T} {{\frac{{dxdy}}{{\ln xy}}}}} .
\]

\noindent
Substituting the series (\ref{eq13}) into the first integral, we integrate termwise
and obtain the series in the desired formula. Letting $y = {{u}
\mathord{\left/ {\vphantom {{u} {x}}} \right. \kern-\nulldelimiterspace}
{x}}$ in the second integral gives

\[
{\int\!\!\!\int_{T} {{\frac{{dxdy}}{{\ln xy}}}}}  = {\int_{0}^{1}
{\,{\frac{{1}}{{x}}}{\int_{x - x^{2}}^{x} {{\frac{{du}}{{\ln u}}}dx}}} }  =
{\int_{0}^{1} {\,{\frac{{\li(x) - \li(x - x^{2})}}{{x}}}\,}} dx,
\]
 
\noindent
and the following calculation (see [\textbf{8}, Section 6.212]) completes
the proof:

\[
{\int_{0}^{1} {\,{\frac{{\li(x)}}{{x}}}\,}} dx = {\mathop {\lim} \limits_{q
\to 0}} {\int_{0}^{1} {\,{\frac{{\li(x)}}{{x^{q + 1}}}}\,}} dx = {\mathop
{\lim} \limits_{q \to 0}} {\frac{{\ln (1 - q)}}{{q}}} = - 1.
\tag*{$\Box$}
\]

\noindent
\textbf{Proposition 2} \textit{The following formula for Euler's constant is valid}:

\[
\gamma = {\int_{0}^{\infty}  {{\sum\limits_{k = 2}^{\infty}
{{\frac{{1}}{{k^{2} \binom{t + k}{k}}}dt}}}} } .
\]

\noindent
\textit{Proof.} In (\ref{eq17}) we replace \textit{xy} with $1 - y$. Substituting the result into (\ref{eq9}), we
change the order of integration and get
\begin{align*}
 \gamma & = {\int_{0}^{\infty}  {{\int_{0}^{1} {{\int_{1 - y}^{1} {{\frac{{1 -
x}}{{xy}}}(1 - y)^{t}dxdy}}} } dt}}  = {\int_{0}^{\infty}  {{\int_{0}^{1}
{{\frac{{ - \ln (1 - y) - y}}{{y}}}(1 - y)^{t}}} dydt}}  \\
 \\
&  = {\int_{0}^{\infty}  {{\sum\limits_{k = 2}^{\infty}
{{\frac{{1}}{{k}}}{\int_{0}^{1} {y^{k - 1}(1 - y)^{t}dydt}}} }} } . 
\end{align*}
Using the proof of Proposition 1, we obtain the desired formula. $\qed$

\section{Integrals over higher-dimensional analogs of $T$}

There are several ways to generalize the triangle $T$ and the integral $I_{n}
$. The simplest generalization of $T$ is the polytope

\[
V_{m} : = {\left\{ {{\left. {(x_{1} ,x_{2} ,\dots, x_{m} ) \in
[0,1]^{m}} \right|}x_{1} + x_{j} \ge 1,\;j = 2, \dots, m} \right\}}.
\]

\noindent
\textbf{Theorem 6} \textit{For} $m \ge 2$\textit{ and} $n \ge 0$
\textit{, the integral}

\[
K_{m,n} : = \int { \cdots {\int_{V_{m}}  {{\frac{{( - \ln (x_{1}
x_{2} \cdots x_{m} ))^{n}}}{{x_{1} x_{2} \cdots x_{m}
}}}}}}  dx_{1} dx_{2} \cdots dx_{m}
\]

\noindent
\textit{is equal to an integer linear combination of multiple zeta values
of weight} $m + n$ \textit{, namely,}

\begin{equation}
\label{eq18}
K_{m,n} = n!{\sum\limits_{
 {k_{1} \ge 0,\dots,k_{m} \ge 0} \atop
 {k_{1} + \cdots + k_{m} = n} 
 } {{\frac{{(k_{2} + \cdots + k_{m} + m -
1)!}}{{(k_{2} + 1)! \cdots (k_{m} + 1)!}}}}} \,\zeta (k_{2} +
\cdots + k_{m} + m,\{1\}_{k_{1}}  ).
\end{equation}

\noindent
\textit{It is also equal to a polynomial of several variables with rational coefficients in values of the Riemann zeta function at integers.}

\bigskip

\noindent
\textit{Proof.} Expanding $( - \ln (x_{1} x_{2} \cdots x_{m} ))^{n} = ( - \ln
x_{1} - \ln x_{2} - \cdots - \ln x_{m} )^{n}$ gives
\begin{align*}
K_{m,n} & = {\sum\limits_{
{k_{1} \ge 0,\dots,k_{m} \ge 0} \atop
 {k_{1} + \cdots + k_{m} = n}
} {\frac{{n!}}{{k_{1}! k_{2}! \cdots k_{m}!}}}} \\ 
& \times
{\int_{0}^{1} {{\frac{{( - \ln x_{1} )^{k_{1}}} }{{x_{1}}} }  
\left(
{{\int_{1 - x_{1}} ^{1} {{\frac{{( - \ln x_{2} )^{k_{2}}} }{{x_{2}
}}}\,dx_{2}}}
\cdots {\int_{1 - x_{1}} ^{1} {{\frac{{( - \ln
x_{m} )^{k_{m}}} }{{x_{m}}} }\,dx_{m}}} }  \right)\,dx_{1}}}.
\end{align*}

\noindent
Since
\[
{\int_{1 - x_{1}} ^{1} {{\frac{{( - \ln x)^{k}}}{{x}}}\,dx}}  = {\frac{{( -
\ln (1 - x_{1} ))^{k + 1}}}{{k + 1}}},
\]

\noindent
we get
\[
K_{m,n} = {\sum\limits_{
{k_{1} \ge 0, \dots, k_{m} \ge 0} \atop
 {k_{1} + \cdots + k_{m} = n}
} {{\frac{{n!}}{{k_{1} !(k_{2} + 1)! \cdots (k_{m} +
1)!}}}}} {\int_{0}^{1} {{\frac{{( - \ln x_{1} )^{k_{1}}} }{{x_{1}}} }( - \ln
(1 - x_{1} ))^{k_{2} + \cdots + k_{m} + m - 1}\,dx_{1}}}  {\rm
.}
\]

\noindent
Evaluating the last integral using formulas (\ref{eq4}) and (\ref{eq12}), the theorem
follows. $\qed$

\bigskip

Taking $m = 2$, the polytope $V_{2} $ is the triangle $T$, the integral
$K_{2,n} $ is the triangle integral $I_{n} $, and Theorem 6 reduces to
Theorem 1. In particular, formula (\ref{eq18}) for $K_{m,n} $ is a weighted version
of formula (\ref{eq2}) for $I_{n} $.

\bigskip

\noindent
\textbf{Corollary 6} \textit{If} $m \ge 2$\textit{, then} $K_{m,0} = (m - 1)!\zeta (m).$

\bigskip

\noindent
\textit{Proof.} Taking $n = 0$ in Theorem 6 forces $k_{1} = k_{2} = \cdots =
k_{m} = 0$ in (\ref{eq18}). $\qed$

\bigskip

There is another, more natural proof of Corollary 6, one that does not use
Theorem 6.

\bigskip

\noindent
\textit{Second proof of Corollary} 6. We use the representation

\[
\zeta (m) = \int { \cdots {\int_{[0,1]^{m}} {{\frac{{dx_{1} \cdot
\cdot \cdot dx_{m}}} {{1 - x_{1} \cdots x_{m}}} }}}}  .
\]

\noindent
(To prove this formula, expand the integrand in a geometric series and
integrate termwise.) We perform the change of variables

\[
x_{1} = y_{m} ,\;\;x_{2} = {\frac{{y_{m - 1}}} {{y_{m}}} },\;\;x_{3} =
{\frac{{y_{m - 2}}} {{y_{m - 1}}} },\dots\;,x_{m - 1} = {\frac{{y_{2}
}}{{y_{3}}} },\;\;x_{m} = {\frac{{1 - y_{1}}} {{y_{2}}} }.
\]

\noindent
(For $m = 2$, compare this with the transformation of $I_{0} $ into Beukers'
integral for $\zeta (2)$ in the Introduction.) We get

\[
\zeta (m) = K'_m := \int { \cdots \int {{\frac{{dy_{1} \cdots
dy_{m}}} {{y_{1} \cdots y_{m}}} }}}  ,
\]

\noindent
where the integral is over the polytope defined by $1 \ge y_{m} \ge y_{m -
1} \ge \cdots \ge y_{2} \ge 0$, $y_{2} + y_{1} \ge 1$, $y_{1} \le
1$. By symmetry, we may interchange the variables $y_{i}$ and $y_{j}$ if
$i > j > 1$. Using all permutations of $y_{2}, \dots, y_{m}$, we
arrive at

\[
(m - 1)!\zeta (m) = \int { \cdots {\int_{V_{m}^{'}}
{{\frac{{dy_{1} \cdots dy_{m}}} {{y_{1} \cdots y_{m}
}}}}}}  ,
\]

\noindent
where the integral is over
\[
V_{m}^{'} : = {\left\{ {{\left. {(y_{1} ,y_{2} ,\dots\ y_{m} ) \in
[0,1]^{m}} \right|}y_{1} + y_{j} \ge 1,\;j = 2, \dots, m} \right\}}.
\]
This is the integral $K_{m,0} $, and the second proof is complete. $\qed$

\bigskip

\noindent
\textbf{Remark 3} The integral $K'_m$ is equivalent to the
Chen (Drinfeld-Kontsevich) integral [\textbf{16}, Section 1]
$$
\int \cdots \int \frac {dY_1 dY_2 \cdots dY_m} {(1-Y_1) Y_2 \cdots Y_m}
$$
over the polytope $1 \ge Y_m \ge \cdots \ge Y_1 \ge 0$: set
$y_1 = 1 - Y_1$ and $y_j = Y_j$ for $j=2,\ldots,m$.

\bigskip

If we rewrite (\ref{eq18}) as

\[
K_{m,n} = n!{\sum\limits_{p = 0}^{n} {a_{m,p} \zeta (m + p,\{1\}_{n - p} )}
},
\]

\noindent
where $a_{m,p} $ denotes the sum of the multinomial coefficients

\[
a_{m,p} : = {\sum\limits_{
 {k_{2} \ge 0,\dots,k_{m} \ge 0} \atop
 {k_{2} + \cdots + k_{m} = p}
} {{\frac{{(k_{2} + \cdots + k_{m} + m -
1)!}}{{(k_{2} + 1)! \cdots (k_{m} + 1)!}}}}} ,
\]

\noindent
then for any fixed $m \ge 2$ one can derive a closed expression for $a_{m,p}
$. For example,

\[
a_{2,p} = 1,
\quad
a_{3,p} = 4 \cdot 2^{p} - 2,
\quad
a_{4,p} = 27 \cdot 3^{p} - 24 \cdot 2^{p} + 3.
\]

\noindent
In general, we have

\bigskip

\noindent
\textbf{Proposition 3} \textit{Fix} $m \ge 2$\textit{ and} $p \ge 0$
\textit{. Then the integers} $a_{m,p}$, $a_{m - 1,p + 1}$
,\dots, $a_{2,p + m - 2}$ \textit{satisfy the recurrence}

\[
{\sum\limits_{t = 0}^{m - 2} {\binom{m - 1}{t}
a_{m - t,p + t}}}   = (m - 1)^{m + p - 1}.
\]

\noindent
\textit{Proof.} First note that if we denote

\[
S_{m,p} : = {\left\{ \left. (k_{2}, \dots, k_{m}) \in {\Z^{m-1}}
\right| k_{j} \ge - 1, \; j = 2, \dots, m; \; k_{2} + \cdots
+ k_{m} = p \right\}},
\]

\noindent
then

\[
{\sum\limits_{S_{m,p}}  {{\frac{{(k_{2} + \cdots + k_{m} + m -
1)!}}{{(k_{2} + 1)! \cdots (k_{m} + 1)!}}}}}  =
{\sum\limits_{
 {l_{2} \ge 0,\dots,l_{m} \ge 0} \atop
 {l_{2} + \cdots + l_{m} = p + m - 1}
} {{\frac{{(l_{2} + \cdots + l_{m} )!}}{{l_{2} !
\cdots l_{m} !}}}}}  = (m - 1)^{p + m - 1}.
\]

\noindent
Now note that if $S_{m,p,t} $ is the subset of $S_{m,p}$ consisting of
those $(m - 1)$-tuples $(k_{2}, \dots, k_{m})$ with exactly $t$ numbers
among the $k_{j} $ equal to $ - 1$, then

\begin{align*}
\sum\limits_{S_{m,p,t}}
\frac{(k_{2} + \cdots + k_{m} + m - 1)!}
{(k_{2} + 1)! \cdots (k_{m} + 1)!} & = 
\binom{m - 1}{t}
\sum\limits_{
 {l_{2} \ge 0,\dots,l_{m} \ge 0} \atop
 {l_{2} + \cdots + l_{m - t} = p + t} 
}
\frac{(l_{2} + \cdots + l_{m} + m - 1 - t)!}
{(l_{2} + 1)! \cdots (l_{m} + 1)!} \\
& = \binom{m - 1}{t} a_{m - t,p + t}.
\end{align*}

\noindent
Finally, since $S_{m,p} $ is the disjoint union

\[
S_{m,p} = {\bigcup\limits_{t = 0}^{m - 2} {S_{m,p,t}}}  ,
\]

\noindent
the proposition follows. $\qed$

\bigskip

Another generalization of the triangle $T$ is the polytope

\[
W_{m} : = {\left\{ {{\left. {(x_{1}, \dots, x_{m} ) \in [0,1]^{m}}
\right|}x_{i} + x_{j} \ge 1,\;1 \le i < j \le m} \right\}}.
\]

\noindent
Note that it is symmetric in all variables, unlike $V_{m} $.

We first generalize the triangle integral $I_{0} $ to an integral over
$W_{m} $. Then we extend to a generalization of $I_{n} $ over $W_{m} $ for
all $n \ge 0$.

Recall that, for all complex $s$ and all $z$ with $|z| < 1$,
the \textit{polylogarithm} $\Li_{s} (z)$ is defined by the convergent series

\[
\Li_{s} (z): = {\sum\limits_{r = 1}^{\infty}  {{\frac{{z^{r}}}{{r^{s}}}}}
}.
\]

\bigskip

\noindent
\textbf{Theorem 7} \textit{If} $m \ge 2$\textit{, then the integral}
\[
L_{m} : = \int { \cdots {\int_{W_{m}}  {{\frac{{dx_{1} \cdot
\cdot \cdot dx_{m}}} {{x_{1} \cdots x_{m}}} }}}}
\]
\textit{is equal to a polynomial of several variables with integer coefficients
in the values} $\ln 2$, $\zeta (m)$,
\textit{and} $\Li_{s} (1/2)$
\textit{ for} $s = 2, 3, \dots, m$ \textit{, namely,}
\[
L_{m} = m!\zeta (m) - (m - 1)\ln ^{m}2 - m!{\sum\limits_{p = 0}^{m - 2}
{{\frac{{\ln ^{p}2}}{{p!}}}}} \Li_{m - p} \left( {{\frac{{1}}{{2}}}}
\right).
\]
\noindent
\textit{Proof.} The symmetry of $W_{m} $ yields

\[
L_{m} = m!\int { \cdots {\int_{W_{m}^{'}}  {{\frac{{dx_{1} \cdot
\cdot \cdot dx_{m}}} {{x_{1} \cdots x_{m}}} }}}}  ,
\]

\noindent
where $W_{m}^{'} $ is the polytope defined by $0 \le x_{1} \le x_{2} \le
\cdots \le x_{m} \le 1$ and $x_{1} + x_{2} \ge 1$. Integrating
consecutively with respect to $x_{m}$, $x_{m - 1}$, \dots, $x_{2}$,
and setting $x = x_{1}$, $y = x_{2}$, we obtain

\[
L_{m} = m(m - 1){\int\!\!\!\int_{H} {{\frac{{( - \ln y)^{m - 2}}}{{xy}}}}
}dxdy,
\]

\noindent
where $H$ is the triangle defined by $0 \le x \le y \le 1$ and $x + y \ge 1$.
(Thus $H$ is the upper half of the triangle $T$ when bisected by the line $y =
x$.) Since $H$ is also defined by ${{1} \mathord{\left/ {\vphantom {{1} {2}}}
\right. \kern-\nulldelimiterspace} {2}} \le y \le 1$ and $1 - y \le x \le
y$, we see that

\begin{align*}
 L_{m} & = m(m - 1){\int_{{{1} \mathord{\left/ {\vphantom {{1} {2}}} \right.
\kern-\nulldelimiterspace} {2}}}^{1} {{\frac{{( - \ln y)^{m - 2}(\ln y - \ln
(1 - y))}}{{y}}}\,dy}}  \\
& = (m - 1)\left( { - (\ln 2)^{m} + m{\int_{{{1} \mathord{\left/ {\vphantom
{{1} {2}}} \right. \kern-\nulldelimiterspace} {2}}}^{1} {{\frac{{( - \ln
y)^{m - 2}( - \ln (1 - y))}}{{y}}}\,dy}}}  \right).
\end{align*}

\noindent
The series expansion

\[
{\frac{{ - \ln (1 - y)}}{{y}}} = {\sum\limits_{r = 1}^{\infty}
{{\frac{{y^{r - 1}}}{{r}}}}}
\]

\noindent
yields

\[
{\int_{{{1} \mathord{\left/ {\vphantom {{1} {2}}} \right.
\kern-\nulldelimiterspace} {2}}}^{1} {{\frac{{( - \ln y)^{m - 2}( - \ln (1 -
y))}}{{y}}}\,dy}}  = {\sum\limits_{r = 1}^{\infty}  {{\frac{{1}}{{r}}}}
}{\int_{{{1} \mathord{\left/ {\vphantom {{1} {2}}} \right.
\kern-\nulldelimiterspace} {2}}}^{1} {y^{r - 1}( - \ln y)^{m - 2}\,dy}} {\rm
.}
\]

\noindent
Now consider the identity

\[
{\int_{{{1} \mathord{\left/ {\vphantom {{1} {2}}} \right.
\kern-\nulldelimiterspace} {2}}}^{1} {y^{t + r - 1}\,dy}}  = {\frac{{1}}{{t
+ r}}} - {\frac{{1}}{{(t + r)2^{t + r}}}}.
\]

\noindent
If we differentiate $l$ times with respect to $t$, and then set $t = 0$, we obtain

\[
{\int_{{{1} \mathord{\left/ {\vphantom {{1} {2}}} \right.
\kern-\nulldelimiterspace} {2}}}^{1} {y^{r - 1}( - \ln y)^{l}\,dy}}  =
{\frac{{l!}}{{r^{l + 1}}}} - {\frac{{l!}}{{2^{r}}}}{\sum\limits_{p =
0}^{l} {{\frac{{\ln ^{p}2}}{{r^{l - p + 1}p!}}}}} .
\]

\noindent
Putting $l = m - 2$, the theorem follows. $\qed$

\bigskip

\noindent
\textbf{Example 1} Taking $m = 2$ gives

\[
L_{2} = {\int\!\!\!\int_{W_{2}}  {{\frac{{1}}{{xy}}}}}  = 2\zeta (2) - \ln
^{2}2 - 2\Li_{2} \left( {{\frac{{1}}{{2}}}} \right).
\]

\noindent
On the other hand, $L_{2} = I_{0} = \zeta (2)$. This proves Euler's formula
for the dilogarithm at ${{1} \mathord{\left/ {\vphantom {{1} {2}}} \right.
\kern-\nulldelimiterspace} {2}}$ [\textbf{5}, Section 1.2], [\textbf{7}, pp.
43-45], [\textbf{10}, Section 1.4]:

\begin{equation}
\label{eq19}
\Li_{2} \left( {{\frac{{1}}{{2}}}} \right) = {\sum\limits_{r = 1}^{\infty}
{{\frac{{1}}{{r^{2}2^{r}}}}}}  = {\frac{{\zeta (2)}}{{2}}} - {\frac{{\ln
^{2}2}}{{2}}}.
\end{equation}

Now take $m = 3$. Using Landen's formula for the trilogarithm at ${{1}
\mathord{\left/ {\vphantom {{1} {2}}} \right. \kern-\nulldelimiterspace}
{2}}$ [\textbf{10}, Equation 6.12],

\begin{equation}
\label{eq20}
\Li_{3} \left( {{\frac{{1}}{{2}}}} \right) = {\frac{{7\zeta (3)}}{{8}}} -
{\frac{{\pi ^{2}\ln 2}}{{12}}} + {\frac{{\ln ^{3}2}}{{6}}}.
\end{equation}
 
\noindent
we get

\[
L_{3} = {\int\!\!\!\int\!\!\!\int_{W_{3}}  {{\frac{{dx_{1} dx_{2} dx_{3}
}}{{x_{1} x_{2} x_{3}}} }}}  = 6\zeta (3) - 2\ln ^{3}2 - 6\Li_{3} \left(
{{\frac{{1}}{{2}}}} \right) - 6\ln 2\Li_{2} \left( {{\frac{{1}}{{2}}}}
\right) = {\frac{{3}}{{4}}}\zeta (3).
\]

\noindent
Thus, surprisingly, $L_{3} $ and $I_{1} $ are both rational multiples of
$\zeta (3)$.

Finally, setting $m = 4$ and using the formulas for $\Li_{2} ({{1}
\mathord{\left/ {\vphantom {{1} {2}}} \right. \kern-\nulldelimiterspace}
{2}})$ and $\Li_{3} ({{1} \mathord{\left/ {\vphantom {{1} {2}}} \right.
\kern-\nulldelimiterspace} {2}})$, we obtain

\[
L_{4} = {\frac{{4}}{{15}}}\pi ^{4} - \ln ^{4}2 + \pi ^{2}\ln ^{2}2 - 21\zeta
(3)\ln 2 - 24\Li_{4} \left( {{\frac{{1}}{{2}}}} \right).
\]

\bigskip

We now generalize $I_{n} $ to an integral over the polytope $W_{m} $. First,
we extend the definition of the polylogarithm $\Li_{s} (z)$ by defining the
\textit{multiple polylogarithm}

\[
\Li_{s_{1}, \dots, s_{l}}  (z): = {\sum\limits_{n_{1} > n_{2} > \cdot \cdot
\cdot > n_{l} > 0} {{\frac{{z^{n_{1}}} }{{n_{1}^{s_{1}}  \cdots
n_{l}^{s_{l}}} } }}} .
\]

\noindent
\textbf{Theorem 8} \textit{If} $m \ge 2$\textit{ and} $n \ge 0$
\textit{, then the integral}

\[
M_{m,n} : = \int { \cdots {\int_{W_{m}}  {{\frac{{( - \ln (x_{1}
\cdots x_{m} ))^{n}}}{{x_{1} \cdots x_{m}}} }dx_{1}
\cdots dx_{m}}} }
\]

\noindent
\textit{is equal to a polynomial of several variables with rational coefficients
in the values} $\ln 2$, $\zeta (a,\{1\}_{m + n - a} )$
\textit{ with} $m \le a \le m + n$\textit{, and} $\Li_{b,\{1\}_{c}}
({{1} \mathord{\left/ {\vphantom {{1} {2}}} \right.
\kern-\nulldelimiterspace} {2}})
$\textit{ with} $b + c \le m + n,\;b \ge 2,\;0 \le c \le n. $

\noindent
\textit{Explicitly, if} $A(k_{2} ): = {\frac{{1}}{{k_{2} !}}}$\textit{ and if}
\[
A(k_{2} , \dots, k_{m} ): = {\frac{{1}}{{k_{2} ! \cdots k_{m}
!}}} \cdot {\frac{{1}}{{(k_{m} + 1)(k_{m - 1} + k_{m} + 2) \cdots
(k_{3} + \cdots + k_{m} + m - 2)}}}
\]

\noindent
\textit{for} $m \ge 3$\textit{, then}
\begin{align*}
M_{m,n} = m!n!{\sum\limits_{
 {k_{1} \ge 0,\dots,k_{m} \ge 0} \atop
 {k_{1} + \cdots + k_{m} = n} 
} {A(k_{2}, \dots, k_{m} ){ \biggl[ {(k_{2} + \cdots +
k_{m} + m - 2)!\,\zeta (k_{2} + \cdots + k_{m} + m,\{1\}_{k_{1}}
)} }}} \\
 - {\frac{{\ln ^{m + n}2}}{{(k_{1} + 1)!(m + n)}}} - (k_{2} + \cdot \cdot
\cdot + k_{m} + m - 2)!{\sum\limits_{p = 0}^{k_{2} + \cdots +
k_{m} + m - 2} {{\frac{{\ln ^{p}2}}{{p!}}}}} { {\Li_{k_{2} + \cdot
\cdot \cdot + k_{m} + m - p,\{1\}_{k_{1}}}   \left( {{{1} \mathord{\left/
{\vphantom {{1} {2}}} \right. \kern-\nulldelimiterspace} {2}}} \right)}
\biggr]}.
\end{align*}

\noindent
\textit{Proof.} The proof is similar to that of Theorem 7 (the case $n = 0)$.
$\qed$

\bigskip

Notice the equality of the integrals $M_{2,n} = I_{n} $.

As an application of Theorem 8, we obtain the following relation between
certain multiple polylogarithm values and multiple zeta values. (The
relation can also be deduced from the H\"{o}lder convolution formula in
[\textbf{4}, Equation (7.2)].)

\bigskip

\noindent
\textbf{Corollary 7} \textit{If} $n \ge 0$\textit{, then}
\[
{\sum\limits_{k = 0}^{n} {{\sum\limits_{p = 0}^{n - k} {{\frac{{\ln
^{p}2}}{{p!}}}}} \Li_{n - k + 2 - p,\{1\}_{k}}  \left( {{\frac{{1}}{{2}}}}
\right)}}  = {\frac{{1 - 2^{n + 1}}}{{(n + 2)!}}}\ln ^{n + 2}2 +
{\frac{{1}}{{2}}}{\sum\limits_{k = 0}^{n} {\zeta (n - k + 2,\{1\}_{k} )}
}.
\]

\noindent
\textit{Proof.} In Theorem 8, take $m = 2$ and set $k_{1} = k$, so that $k_{2} = n - k$.
Then $A(k_{2} ) = {\frac{{1}}{{(n - k)!}}}$ and
\[
M_{2,n} = 2n!{\sum\limits_{k = 0}^{n} {{\left[ {\zeta (n - k + 2,\{1\}_{k} )
- {\frac{{\ln ^{n + 2}2}}{{(k + 1)!(n - k)!(n + 2)}}} - {\sum\limits_{p =
0}^{n - k} {{\frac{{\ln ^{p}2}}{{p!}}}}} \Li_{n - k + 2 - p,\{1\}_{k}}
\left( {{\frac{{1}}{{2}}}} \right)} \right]}}} .
\]

\noindent
Now substitute $M_{2,n} = I_{n} $, and apply Theorem 1 (i) and the identity

\[
{\sum\limits_{k = 0}^{n} {{\frac{{1}}{{(k + 1)!(n - k)!}}}}}  = {\frac{{2^{n
+ 1} - 1}}{{(n + 1)!}}}.
\tag*{$\Box$}
\]

\noindent
\textbf{Example 2} The case $n = 0$ is Euler's formula (\ref{eq19}) for $\Li_{2}
\left( {{{1} \mathord{\left/ {\vphantom {{1} {2}}} \right.
\kern-\nulldelimiterspace} {2}}} \right)$. Taking $n = 1$ and substituting
$\zeta (2,1) = \zeta (3)$ gives the relation

\[
\Li_{3} \left( {{\frac{{1}}{{2}}}} \right) + \Li_{2} \left(
{{\frac{{1}}{{2}}}} \right)\ln 2 + \Li_{2,1} \left( {{\frac{{1}}{{2}}}}
\right) = - {\frac{{\ln ^{3}2}}{{2}}} + \zeta (3),
\]

\noindent
which is a special case of [\textbf{4}, Equation (7.3)]. Using the values of
$\Li_{2} ({{1} \mathord{\left/ {\vphantom {{1} {2}}} \right.
\kern-\nulldelimiterspace} {2}})$ and $\Li_{3} ({{1} \mathord{\left/
{\vphantom {{1} {2}}} \right. \kern-\nulldelimiterspace} {2}})$ in (\ref{eq19}) and
(\ref{eq20}), we get the formula

\[
\Li_{2,1} \left( {{\frac{{1}}{{2}}}} \right) = {\sum\limits_{r = 2}^{\infty}
{{\frac{{1}}{{r^{2}2^{r}}}}\left( {1 + {\frac{{1}}{{2}}} + \cdots
+ {\frac{{1}}{{r - 1}}}} \right)}}  = {\frac{{\zeta (3)}}{{8}}} -
{\frac{{\ln ^{3}2}}{{6}}}.
\]

\noindent
Finally, adding $\Li_{3} ({{1} \mathord{\left/ {\vphantom {{1} {2}}} \right.
\kern-\nulldelimiterspace} {2}})$ recovers Ramanujan's summation
[\textbf{1}, p. 258]

\[
{\sum\limits_{r = 1}^{\infty}  {{\frac{{1}}{{r^{2}2^{r}}}}\left( {1 +
{\frac{{1}}{{2}}} + \cdots + {\frac{{1}}{{r}}}} \right)}}  =
\zeta (3) - {\frac{{\pi ^{2}\ln 2}}{{12}}}.
\]

\section{References}

\noindent
1. B. C. Berndt, \textit{Ramanujan's Notebooks}, Part I, Springer-Verlag, New York, 1985.

\noindent
2. F. Beukers, A note on the irrationality of $\zeta (2)$ and $\zeta (3)$,
\textit{Bull. London Math. Soc.} \textbf{11} (1979) 268-272.

\noindent
3. J. M. Borwein, D. M. Bradley, and D. J. Broadhurst, Evaluations of $k$-fold
Euler/Zagier sums: a compendium of results for arbitrary $k$, \textit{Electron. J. Combin.} \textbf{4}
(1997) no. 2, R5.

\noindent
4. J. M. Borwein, D. M. Bradley, D. J. Broadhurst, and P. Lisonek, Special
values of multiple polylogarithms, \textit{Trans. Amer. Math. Soc.} \textbf{353} (2001) 907-941.

\noindent
5. P. Cartier, Fonctions polylogarithmes, nombres polyz\^{e}tas et groupes
pro-unipotents, \textit{Ast\'{e}risque} \textbf{282} (2002) 137-173.

\noindent
6. N. G. de Bruijn, \textit{Asymptotic Methods in Analysis}, Dover Publications, New York, 1981.

\noindent
7. W. Dunham, \textit{Euler: The Master of Us All}, Dolciani Mathematical Expositions No. 22, Mathematical
Association of America, Washington, D.C., 1999.

\noindent
8. S. Gradshteyn and I. M. Ryzhik, \textit{Table of Integrals, Series, and Products}, 6th edition, A. Jeffrey and D. Zwillinger,
editors, Academic Press, San Diego, 2000.

\noindent
9. K. S. K\"{o}lbig, J. A. Mignaco, and E. Remiddi, On Nielsen's generalized
polylogarithms and their numerical calculation, \textit{BIT} \textbf{10} (1970) 38-74.

\noindent
10. L. Lewin, \textit{Polylogarithms and Associated Functions}, North-Holland, New York, 1981.

\noindent
11. J. Ser, Sur une expression de la fonction $\zeta (s)$ de Riemann, \textit{C. R. Acad. Sci. Paris S\'{e}r. I Math.}
\textbf{182} (1926) 1075-1077.

\noindent
12. J. Sondow, Criteria for irrationality of Euler's constant, \textit{Proc. Amer. Math. Soc.} \textbf{131}
(2003) 3335-3344.

\noindent
13. J. Sondow, An infinite product for $e^{\gamma} $ via hypergeometric
formulas for Euler's constant, $\gamma $ (2003, preprint); available at
http://arXiv.org/abs/math/0306008.

\noindent
14. J. Sondow, Double integrals for Euler's constant and $\ln
{\frac{{4}}{{\pi}} }$ and an analog of Hadjicostas's formula, \textit{Amer. Math. Monthly} \textbf{112}
(2005) 61-65.

\noindent
15. J. Sondow, A faster product for $\pi $ and a new integral for $\ln
{\frac{{\pi}} {{2}}}$, \textit{Amer. Math. Monthly} \textbf{112} (2005) 729-734.

\noindent
16. M. Waldschmidt, Multiple polylogarithms: an introduction,
in \textit{Number Theory and Discrete Mathematics,} Proceedings of the
International Conference in Honour of Srinivasa Ramanujan,
Chandigarh, 2000, A. K. Agarwal et al., editors,
Birkh\"{a}user, Basel, 2002, pp. 1-12.

\noindent
17. S. Zlobin, On a certain integral over a triangle (2005, preprint); available
at \\ http://arXiv.org/abs/math/0511239.


\end{document}